Since the group G is regular when p > 3 (this follows from the fact that the nilpotency class of G is less than p), it is easy to see from the preceding relations that the p-th power of any element of G lies in Z(G). Thus, when p > 3 the p-group does not satisfy conditions 1) and 2).

Obviously, when p = 3 any group of maximal class having order $3^5$ satisfies 1) and 2) (here as a generator of a cyclic subgroup of index $3^3$ we can choose $S_1$). Blackburn described all nonisomorphic types of groups of maximal class having order $3^5$ (see the theorem, 11-16).

Thus, we have proved the following lemma.

LEMMA 5. Suppose a p-group G of maximal class has order $p^5$. If G satisfies condition 1) and 2), then it is a 3-group of one of the types 11-16 (see the theorem).

COVERINGS OF VERTICAL SEGMENTS UNDER A CONFORMAL MAPPING

V. N. Dubinin

1. Let f(z) be an arbitrary function, regular in the disk $U = \{z: |z| < 1\}$ with expansion $f(z) = a_0 + a_1 z + \ldots$.

The following problem is posed in [1]: does the inequality $|a_1| \leqslant 2l/\pi$ follow from the assumption that f(U) does not cover any vertical segment of length greater than $l$?

In [2] on the basis of the method of symmetrization and the extremal metric principle, an affirmative answer was given to the question posed and its proof was outlined. In the present note we give a detailed proof for schlicht functions, formulating the basic result as a theorem on coverings of vertical segments under a conformal mapping.

We denote by S the class of regular and schlicht functions w = f(z) in the disk U, normalized by the conditions f(0) = 0, f'(0) = 1.

THEOREM. Let f be an arbitrary function of the class S. Then:

1. f(U) covers some vertical segment of length greater than $\pi/2$, except for the case when

$$f(z) = \frac{1}{2} \ln \frac{1+z}{1-z}.$$

2. If the indicated segment belongs to the line Re w = $u_0$ (w = u + iv), then there exists a continuous curve $L = \{w: w = u + iv(u), \ 0 \leqslant u \leqslant u_0 \ (u_0 \leqslant u \leqslant 0), \ v(0) = 0\}$ such that $L \subset f(U)$.

COROLLARY. If the function $f(z) = a_0 + a_1 z + \ldots$ is regular and schlicht in the disk U, and f(U) does not cover any vertical segment of length greater than $l$, then $|a_1| \leqslant 2l/\pi$.

The equality sign holds only for functions of the form

$$f(z) = a_0 + \frac{l}{\pi} \ln \frac{1+\varepsilon z}{1-\varepsilon z}, \quad |\varepsilon| = 1.$$

---

  

In fact, the function $q(z) = \left[f\left(z\frac{|a_1|}{a_1}\right) - a_0\right] / |a_1| \in S$ and by the theorem q(U) contains some vertical segment of length greater than $\pi/2$. But according to the assumption, the length of this segment does not exceed $l/|a_1|$. Hence, $l/|a_1| \geqslant \pi/2$ with the corresponding assertion about the equality sign.

Remark 1. In the case when $f(z) = a_0 + a_1 z + \ldots$ is an arbitrary regular function, defined in the disk U, the problem can be solved by the same scheme as the proof of the theorem, but in place of F(U) it is necessary to consider the Riemann surface of the inverse map.

2. Definition and Auxiliary Results. We introduce the following notation:

$$B(r) = B(r, f) = \{w: w = f(rz),\ z \in U\},\ 0 \leqslant r \leqslant 1,$$
$$\gamma(r) = \gamma(r, f) = \{w: w = f(z),\ |z| = r\},\ 0 \leqslant r < 1,\ w = u + iv.$$

Let $\{\lambda\}$ be some system of lines u = const, containing all tangents to $\gamma(\rho, f)$, $f \in S$, $\rho < 1$. Lines from $\{\lambda\}$ divide the set $B(\rho, f) \cap \{w: u > 0\}$ into a finite number of domains. Those domains which can be joined in the set $B(\rho, f)$ by a curve of the form v = v(u) with the point w = 0, constitute the collection $\{D\}$. Let r(D) be the smallest r, $0 \leqslant r \leqslant \rho$, for which any line u = const, passing through $\overline{D}$ ($D \in \{D\}$), intersects $\overline{B}(r, f)$.

We index the domains of $\{D\}$ in such a way that one has the relations

$$\{w=0\} \in \overline{D}_1,\ r(D_k) = \min_{D \in \{D\} \setminus \cup_{m+1}^{k-1} D_m} r(D),\ k = 2, 3, \ldots, n.$$

Let $\alpha_k$ be the left boundary continuum of the domain $D_k$, lying on a line from $\{\lambda\}$, and let $\beta_k$ be the right one.

Proposition 1. If the line $\beta$: u = u' intersects the set $\overline{D_k \cap B}(r)$, $r < \rho$, then any line $\alpha$: u = u'', lying between $\beta$ and $\alpha_k$, also intersects $\overline{D_k \cap B}(r)$.

Proof. Let $\alpha \cap \overline{D_k \cap B}(r) = \emptyset$. The set $\alpha \cap D_k$ is a cross section of the simply connected domain $B(\rho)$. Hence it divides $B(\rho)$ into two domains, where $\alpha_k \cap B(\rho)$ and $\beta \cap D_k$ lie in different domains. However, $\alpha_k \cap B(\rho)$ can be joined in $B(\rho)$ by a curve of the form v = v(u) with the point w = 0, and then by a continuous curve $\gamma \subset \overline{B}(r)$ with $\beta \cap \overline{D_k \cap B}(r)$, which in view of the assumption made does not intersect $\alpha \cap D_k$. Contradiction.

COROLLARY. Upon adding to $\{\lambda\}$ new lines, the number of sets $D_k$, satisfying the conditions $\overline{D_k \cap B}(r) \cap \alpha_k \neq \emptyset$, $\overline{D_k \cap B}(r) \cap \beta_k = \emptyset$ does not increase.

Proposition 2.

$$\sum_{m=k}^{n} \operatorname{mes} B(r) \cap \beta_m \leqslant \sum_{m=k+1}^{n} \operatorname{mes} B(r) \cap \alpha_m. \tag{1}$$

Proof. Let $\beta_{m_1}, \beta_{m_2}, \ldots, \beta_{m_s}$ be all possible continua with indices $m_{s'} \geqslant k$, lying on the line $u = u_0$; $\alpha_{l_1}, \alpha_{l_2}, \ldots, \alpha_{l_t}$ be continua having nonempty intersection with $\bigcup_{s'=1}^{s} \beta_{m_{s'}} \setminus \gamma(\rho)$. Then

$$\sum_{s'=1}^{s} \operatorname{mes} B(r) \cap \beta_{m_{s'}} \leqslant \sum_{l'=1}^{t} \operatorname{mes} B(r) \cap \alpha_{l_{l'}}. \tag{1'}$$

We note however that if $(\beta_m \cap \alpha_l) \setminus \gamma(\rho) \neq \emptyset$, then $l \geqslant m$. In fact, let $l < m$. Then, on the one hand, $\beta_m$ divides $B(\rho)$ into two domains containing $D_m$ and $D_l$, respectively, and on the other, $D_m$ can be joined by a curve of the form v = v(u) with the point w = 0, and then by a continuous curve $\gamma$, $\gamma \subset B(r(D_l))$ with $D_l$. But by Proposition 1, $B(r(D_l)) \cap \beta_m = \emptyset$. Contradiction.

Summing inequalities (1') over all lines $u = u_0$ from the system $\{\lambda\}$, we get (1).

Everywhere below, G* will denote the result of Steiner symmetrization of the set G with respect to the real axis (see, e.g., [1]). Let $B(r) \cap D_k \neq \emptyset$. The interval $[a_k, b_k(r)] = \{u: w \in \overline{B(r) \cap D_k}\}$ can be divided into a finite number of intervals, in each of which the curve $\gamma(r) \cap D_k$ is defined by equations of the form

$$v = v_{k,c}(u),\ 1 \leqslant c \leqslant 2C,\ C \geqslant 1,$$
$$v_{k,c+1}(u) \leqslant v_{k,c}(u),\ c = 1, 2, \ldots, 2C - 1,$$

and the equation of the "upper" part of the boundary of the set $(B(r) \cap D_k)^*$ has the form $v = v_k^*(u)$.

Let the function $w' = F_k(w)$, $F_k: D_k \mapsto D_k^*$, correspond to each point $w \in \gamma(r) \cap D_k$, $\operatorname{Re} w = u$, $\operatorname{Im} w = v_{k,c}(u)$, $c = 1, 2, \ldots, 2C$, a point $w'$, $\operatorname{Re} w' = u$, $\operatorname{Im} w' = v_k^*(u)$. $H_k(w)$, $k = 1, 2, \ldots, n$, is the translation of the set $\overline{D_k^*}$ along the real axis such that $H_1(a_1) = 0$, $H_{k-1}(b_{k-1}(\rho)) = H_k(a_k)$, $k = 2, 3, \ldots, n$. The curves $v = v_k^*(u)$ under the maps $H_k(w)$ are carried into curves

$$\gamma_k(r): v = v_k(u), \quad a_k' \leqslant u \leqslant b_k'(r),$$
$$H_k(a_k) = a_k', \quad H_k(b_k(r)) = b_k'(r).$$
$$w = p(z) = \frac{1}{2} \ln[(1+z)/(1-z)], \quad p(z) \in S.$$
$$\mu(w) = \frac{1}{2\pi} \left| \frac{d}{dw} \ln p^{-1}(w) \right|, \quad |u| < \infty, \quad |v| \leqslant \pi/4.$$

In what follows we shall assume that the domain B($\rho$) does not contain a vertical segment of length greater than $\pi/2$, and the curves $\gamma(r, f)$ are convex for all $r \leqslant r^*$, $0 < r^*$.

$$\mu^+(w) = \begin{cases} \mu(H_k(F_k(w))), & w \in D_k, \quad k = 1, 2, \ldots, n, \\ 0, & w \in B(\rho) \setminus \bigcup_{k=1}^{n} D_k. \end{cases}$$

By a lift of the curve $\gamma_k(r)$ we shall mean a translation of this curve "upwards" along the axis u = 0 with the simultaneous replacement of the arc of the curve issuing from the line v = $\pi/4$, by the corresponding segments of this line.

<u>Proposition 3.</u> Let $\delta$ be the largest distance between adjacent lines of the system $\{\lambda\}$, and $0 < r_0 < r^*$. Then

$$\int_{\gamma(r)} \mu^+(w) |dw| \geqslant 1/2 - \varepsilon(\delta), \tag{2}$$

where $\varepsilon(\delta)$ is independent of $r$, $r_0 < r < \rho$ and $\lim_{\delta \to 0} \varepsilon(\delta) = 0$.

<u>Proof.</u> Let $\gamma(r) \cap D_k \neq \varnothing$. Applying Minkowski's inequality [3] $2C - 1$ times, we get that for almost all $u \in [a_k, b_k(r)]$

$$\sum_{c=1}^{2C} \sqrt{1 + \left(\frac{\partial v_{k,c}}{\partial u}\right)^2} \sqrt{(2C)^2 + \left(\sum_{c=1}^{2C} (-1)^{c+1} \frac{\partial v_{k,c}}{\partial u}\right)^2} \geqslant 2 \sqrt{1 + \left(\frac{\partial v_k^*}{\partial u}\right)^2}. \tag{3}$$

Equality is possible only in the case $C = 1$, $\frac{\partial v_{k,1}}{\partial u} = \frac{\partial v_k^*}{\partial u} = -\frac{\partial v_{k,2}}{\partial u}$. It follows from (3) that

$$\int_{\gamma(r) \cap D_k} \mu^+(w) |dw| \geqslant 2 \int_{\gamma_k(r)} \mu(w) |dw|.$$

According to Proposition 1, the continuous curve $\gamma_k(r)$ joins the line $u = a_k'$ with the line $u = b_k'(\rho)$ for $k < k'$, the line $u = a_k'$ with the line v = 0 for $k' \leqslant k \leqslant n'$, and $\gamma_k(r) = \varnothing$ for $k > n'$ ($k'$, $n'$, $1 \leqslant k' \leqslant n' \leqslant n$, are certain natural numbers depending on r and $\delta$).

We construct successive lifts of the curves $\gamma_k(r)$, $k = n', n' - 1, \ldots, 1$, from right to left so that we raise the curve $\gamma_k(r)$ by the quantity

$$\sum_{m=k+1}^{n'} v_m(a_m') - \sum_{m=k}^{n'} v_m(b_m'(r)),$$

which is nonnegative by (1). We denote the union of the curves so obtained by $\gamma'(r)$.

Since $\mu(u + iv) > \mu(u + i(v + \Delta v))$, $|u| < \infty$, $0 < v < v + \Delta v < \pi/4$, one has

$$\sum_{k=1}^{n'} \int_{\gamma_k(r)} \mu(w) |dw| \geqslant \int_{\gamma'(r)} \mu(w) |dw|.$$

From the curve $\gamma'(r)$ one gets a continuous curve $\gamma^*(r)$, if its parts $\gamma_k'(r)$, obtained by lifting the curves $\gamma_k(r)$, $k = k' + 1, k' + 2, \ldots, n'$, are translated "leftwards" by certain quantities, not exceeding $\delta(n' - k')$. As a result of such translations, the lengths of the curves $\gamma_k'(r)$ in the metric $\mu(w) |dw|$ do not increase by more than the difference

$$\left| \int_{\gamma_k'(r)} \mu(w) |dw| - \int_{\gamma_k'(r)} \mu(w + \delta(n' - k')) |dw| \right| \leqslant$$

$$\leqslant \max_{w \in B(1, \rho) \setminus B^*(r_0, f)} |\mu(w) - \mu(w + \delta(n' - k'))| \max_{r_0 \leqslant r \leqslant \rho} \int_{\gamma(r)} |dw| = \varepsilon'(\delta).$$

According to the corollary to Proposition 1, $(n' - k') < N(\rho) = \text{const}$ for any $r$, $r_0 \leqslant r \leqslant \rho$, and $\delta$. Also taking account of the uniform continuity of $\mu(w)$ on $B(1, p) \setminus B^*(r_0, f)$ and the boundedness of $|f'(z)|$ in $|z| < \rho$, we get $\lim_{\delta \to 0} \varepsilon(\delta) = 0$, where $\varepsilon(\delta) = 2\varepsilon'(\delta)(n' - k')$. It remains to note that in view of the continuity of the curve $\gamma^*(r)$

$$\int_{\gamma^*(r)} \mu(w)|dw| = \int_{p^{-1}(\gamma^*(r))} (2\pi r)^{-1} (dr^2 + r^2 d\varphi^2)^{1/2} \geqslant \frac{1}{2\pi} \int_0^{\pi/2} d\varphi = 1/4, \quad z = re^{i\varphi}.$$

Proposition 4. Let $B(r_0', p) \subset B^*(r_0, f)$, $r_0' \leqslant r_0 < r^*$. Then

$$\iint_{B(\rho) \setminus B(r_0, f)} (\mu^+(w))^2 du\, dv \leqslant -\frac{1}{4\pi} \ln r_0. \qquad (4)$$

Proof. We denote by $\alpha_k(u)$ the set $(D_k \setminus B(r_0, f)) \cap \{w: \operatorname{Re} w = u\}$. One has the following relation:

$$\iint_{B(\rho) \setminus B(r_0, f)} (\mu^+(w))^2 du\, dv = \sum_{k=1}^n \int_{a_k}^{b_k(\rho)} \int_{\alpha_k(u)} (\mu^+(w))^2 dv\, du = \sum_{k=1}^n \int_{a_k}^{b_k(\rho)} \int_0^\infty \operatorname{mes}\{v: (\mu^+(w))^2 > a,$$

$$w \notin B(r_0, f)\} da\, du = \sum_{k=1}^n \int_{a_k'}^{b_k(\rho)} \int_0^\infty \operatorname{mes}\{v: \mu^2(w) > a,$$

$$w \notin B^*(r_0, f)\} da\, du = \sum_{k=1}^n \int_{H_k(D_k^*) \setminus B^*(r_0, f)} \mu^2(w) du\, dv \leqslant \frac{1}{2} \iint_{B(1, p) \setminus B(r_0', p)} \mu^2(w) du\, dv = -\frac{1}{4\pi} \ln r_0'.$$

Remark 2. If for curves and domains lying in the half plane $\operatorname{Re} w < 0$, by analogy with the arguments given above one defines the function $\mu^-(w)$, then for it (2) and (4) remain valid with the replacement of $\mu^+(w)$ by $\mu^-(w)$.

3. Proof of the Theorem. Let the domain $B(1, f)$, $f \in S$, not contain a vertical segment of length greater than $\pi/2$.

For a given number $r_0$, $0 < r_0 < r^*$, we choose a number $r_0'$ such that $B(r_0', p) \subset B(r_0, f)$ and $\gamma(r_0', p) \cap \gamma(r_0, f) \neq \varnothing$. In this case $B(r_0', p) \subset B^*(r_0, f)$ and $\lim_{r_0 \to 0} r_0'/r_0 = 1$.

Let us assume now that the curve $\gamma^*(r, f) \equiv \partial B^*(r, f)$, $r_0 < r < r^*$, either does not coincide up to a translation along the imaginary axis with the curve $f(r, f)$, or does not coincide with any of the curves $f(r', p)$, $0 < r' < 1$. Then in the first case from (3) and in the second from properties of the logarithmic metric $(2\pi |z|)^{-1}|dz|$ we get

$$\int_{\gamma(r,f)} \mu'(w)|dw| \geqslant 1 + t, \qquad (5)$$

where

$$\mu'(w) = \begin{cases} \mu^+(w), & \operatorname{Re} w \geqslant 0, \\ \mu^-(w), & \operatorname{Re} w < 0, \end{cases}$$

and $t > 0$ is unchanged upon lessening of $\delta$. In view of the continuity of the function $\mu'(w)$ on the set $B(r^*, f) \setminus \{w: u = 0\}$ one can assume that (5) holds for all $r$ from a certain interval $[r_1, r_2]$, $r_0 < r_1 < r_2$.

According to (2) and (5), the metric

$$\mu^*(z)|dz| = \mu'(f(z))|f'(z)||dz|, \quad r_0 < |z| < \rho,$$

satisfies the conditions of Lemma 2.2 [4]. Hence

$$\iint_{r_0 < |z| < \rho} (\mu^*(z))^2 r\, dr\, d\varphi \geqslant \frac{1 - 4\varepsilon(\delta)}{2\pi} \ln \frac{\rho}{r_0} + \frac{t + 2\varepsilon(\delta)}{\pi} \ln \frac{r_2}{r_1},$$

where $z = re^{i\varphi}$, $\lim_{\delta \to 0} \varepsilon(\delta) = 0$.

From Proposition 4 and Remark 2, it follows that

$$-\frac{1}{2\pi} \ln r_0' \geqslant \frac{1 - 4\varepsilon(\delta)}{2\pi} \ln \frac{\rho}{r_0} + \frac{t + 2\varepsilon(\delta)}{\pi} \ln \frac{r_2}{r_1}.$$

Successively letting $\delta$ go to zero, $\rho$ to one, and $r_0$ to zero, we get a contradiction.

To complete the proof of the theorem it remains to consider the case when for some $r < r^*, \gamma(r, f) = \{w: w + c \in \gamma(r', p)\}$, $0 < r' < 1$, $\operatorname{Re} c = 0$. Then $B(r', p) = B^*(r, f)$. Using (2) and (4) and also Lemma 2.2 [4] ($r_0 = r, r_0' = r'$) just as above, and letting δ go to zero, and ρ to one, we get $r' \leqslant r$.

On the other hand, from Schwarz' lemma applied to the composite function $h\{p^{-1}[f(z) +$

$$h(\zeta) = (r')^2 \frac{\zeta - p^{-1}(c)}{(r')^2 + \zeta p^{-1}(c)}, \quad |\zeta| < r',$$

there follows the inequality $r' \geqslant r$, where equality is achieved if and only if $f(z) \equiv p(z)$. Thus, $r = r'$ and $f(z) \equiv p(z)$. The domain $B(1, f)$ in this case is the strip $\{w: |v| < \pi/4\}$. The theorem is proved.

EXACT VALUES OF APPROXIMATION BY HERMITIAN SPLINES OF EVEN DEGREE
IN CLASSES OF DIFFERENTIABLE FUNCTIONS

N. A. Nazarenko and S. V. Pereverzev

1. Denote by $C^r$ ($r = 0, 1, 2, \ldots; C^0 = C$) the linear space of functions f(x), having r continuous derivatives in the interval [0, 1].

Let $\Delta_n = \{0 = x_0 < x_1 < \ldots < x_n = 1\}$ $(n \geqslant 1)$ be a division of interval $[0, 1]$, $h_i = x_i - x_{i-1}$ and $|\Delta_n| = \max\limits_{1 \leqslant i \leqslant n} h_i$. The uniform division, in which $x_i$ = i/n, will be denoted by $\overline{\Delta}_n$.

With every function $f(x) \in C^m$ $(m \geqslant 0)$ we associate functions $S_{2m}(f; x) \in C^m$ and $S_{2m+1}(f; x) \in C^m$, which satisfy the conditions:

a) in any interval $[x_{i-1}, x_i]$ $S_{2m+1}(f; x)$ is an algebraic polynomial of degree 2m + 1, while $S_{2m}(f; x)$ is expressible as

$$s_{2m}(f; x) = \sum_{s=0}^{2m} a_s^{(i)} (x - x_{i-1})^s + a_{2m+1}^{(i)} (x - \bar{x}_{i-1})_+^{2m},$$

where

$$\bar{x}_{i-1} = x_{i-1} + h_i/2, \quad z_+^k = [\max(0, z)]^k;$$

b) $S_{2m+k}^{(j)}(f; x_i) = f^{(j)}(x_i)$ $(k = 0, 1; j = 0, 1, \ldots, m; i = 0, 1, \ldots, n)$.

We know that, for $f(x) \in C^m$, function $S_{2m+k}(f; x)$ (k = 0, 1) exists and is unique (see [1, p. 48] for k = 1; and [2, p. 145], for k = 0). Functions $S_{2m}(f; x)$ and $S_{2m+1}(f; x)$ are called Hermitian splines.

Let

$$e_{2m+k}(f; x) = f(x) - S_{2m+k}(f; x) \quad (k = 0, 1) \tag{1}$$

and let $\mathfrak{M}$ be a class of functions of $C^r$.

We shall consider the determination of the following functions and quantities:

$$e_{2m+k}(\mathfrak{M}; x) = \sup_{f \in \mathfrak{M}} |e_{2m+k}(f; x)| \quad (k = 0, 1), \tag{2}$$